\documentclass[11pt,reqno]{amsart}
\usepackage{tikz-cd}
\usepackage{amsthm}
\usepackage{amssymb}
\usepackage{latexsym}
\usepackage{multicol}
\usepackage{verbatim,enumerate,mathtools}
\usepackage{hyperref}
\usepackage{amsmath, amscd}
\usepackage{soul}
\usepackage{youngtab}
\usepackage{mathrsfs}
\usepackage{tikz-cd}

\advance\textwidth by 1.2in \advance\oddsidemargin by -.6in \advance\evensidemargin by -.6in
\usepackage{tikz}
\usetikzlibrary{decorations.markings}
\tikzstyle{vertex}=[ circle, fill, draw, inner sep=0pt, minimum size=4pt,]
\tikzstyle{edge}= [thick]

\newtheorem*{cor}{Corollary}
\newtheorem*{lem}{Lemma}
\newtheorem*{prop}{Proposition}

\theoremstyle{definition} \newtheorem*{defn}{Definition}
\theoremstyle{definition}
\newtheorem*{thm}{Theorem}
\newtheorem*{thm*}{Theorem}

\newtheorem*{rem}{Remark}

\newenvironment{pf}{\proof}{\endproof}
\newcounter{cnt}
 \makeatletter
\def\mydggeometry{\makeatletter\dg@YGRID=1\dg@XGRID=20\unitlength=0.003pt\makeatother}
\makeatother \theoremstyle{remark}

\numberwithin{equation}{section}

\begin{document}

\newcommand{\Ker}{\operatorname{Ker}}
\newcommand{\thmref}[1]{Theorem~\ref{#1}}
\newcommand{\secref}[1]{Section~\ref{#1}}
\newcommand{\lemref}[1]{Lemma~\ref{#1}}
\newcommand{\propref}[1]{Proposition~\ref{#1}}
\newcommand{\corref}[1]{Corollary~\ref{#1}}
\newcommand{\remref}[1]{Remark~\ref{#1}}
\newcommand{\defref}[1]{Definition~\ref{#1}}
\newcommand{\er}[1]{(\ref{#1})}
\newcommand{\id}{\operatorname{id}}
\newcommand{\ord}{\operatorname{\emph{ord}}}
\newcommand{\sgn}{\operatorname{sgn}}
\newcommand{\wt}{\operatorname{wt}}
\newcommand{\tensor}{\otimes}
\newcommand{\from}{\leftarrow}
\newcommand{\nc}{\newcommand}
\newcommand{\rnc}{\renewcommand}
\newcommand{\dist}{\operatorname{dist}}
\newcommand{\qbinom}[2]{\genfrac[]{0pt}0{#1}{#2}}
\nc{\cal}{\mathcal} \nc{\goth}{\mathfrak} \rnc{\bold}{\mathbb}
\renewcommand{\frak}{\mathfrak}
\newcommand{\supp}{\operatorname{supp}}
\newcommand{\Irr}{\operatorname{Irr}}
\newcommand{\psym}{\mathcal{P}^+_{K,n}}
\newcommand{\psyml}{\mathcal{P}^+_{K,\lambda}}
\newcommand{\psymt}{\mathcal{P}^+_{2,\lambda}}
\renewcommand{\Bbb}{\mathbb}
\nc\bomega{{\mbox{\boldmath $\omega$}}} \nc\bpsi{{\mbox{\boldmath $\Psi$}}}
 
 \nc\bbeta{{\mbox{\boldmath ${\beta_\lambda}$}}}
 \nc\bpi{{\mbox{\boldmath $\pi$}}}
  \nc\bvarpi{{\mbox{\boldmath $\varpi$}}}
  \nc\bnu{{\mbox{\boldmath $\nu$}}}
\nc\bepsilon{{\mbox{\boldmath $\epsilon$}}}

  \nc\bxi{{\mbox{\boldmath $\xi$}}}
\nc\bmu{{\mbox{\boldmath $\mu$}}}
\nc\blambda{{\mbox{\boldmath
$\lambda$}}}

\newcommand{\Tmn}{\bold{T}_{\lambda^1, \lambda^2}^{\nu}}

\newcommand{\lie}[1]{\mathfrak{#1}}
\newcommand{\ol}[1]{\overline{#1}}
\makeatletter
\def\section{\def\@secnumfont{\mdseries}\@startsection{section}{1}%
  \z@{.7\linespacing\@plus\linespacing}{.5\linespacing}%
  {\normalfont\scshape\centering}}
\def\subsection{\def\@secnumfont{\bfseries}\@startsection{subsection}{2}%
  {\parindent}{.5\linespacing\@plus.7\linespacing}{-.5em}%
  {\normalfont\bfseries}}
\makeatother
\def\subl#1{\subsection{}\label{#1}}
 \nc{\Hom}{\operatorname{Hom}}
  \nc{\mode}{\operatorname{mod}}
\nc{\End}{\operatorname{End}} \nc{\wh}[1]{\widehat{#1}} \nc{\Ext}{\operatorname{Ext}}
 \nc{\ch}{\operatorname{ch}} \nc{\ev}{\operatorname{ev}}
\nc{\Ob}{\operatorname{Ob}} \nc{\soc}{\operatorname{soc}} \nc{\rad}{\operatorname{rad}} \nc{\head}{\operatorname{head}}
\def\Im{\operatorname{Im}}
\def\gr{\operatorname{gr}}
\def\mult{\operatorname{mult}}
\def\Max{\operatorname{Max}}
\def\ann{\operatorname{Ann}}
\def\sym{\operatorname{sym}}
\def\loc{\operatorname{loc}}
\def\Res{\operatorname{\br^\lambda_A}}
\def\und{\underline}
\def\Lietg{$A_k(\lie{g})(\bs_\xiigma,r)$}
\def\res{\operatorname{res}}


\title[Generalized Demazure modules and Prime Representations in type $D_n$. ]
{ Generalized Demazure modules and Prime Representations in type $D_n$}

\author[ Chari, Davis, Moruzzi Jr.]{ Vyjayanthi Chari,  Justin Davis and Ryan Moruzzi Jr.}

\address{Department of Mathematics\\ 
  University of California, Riverside\\ 
  900 University Ave., Riverside, CA 92521} \thanks{VC was partially supported by DMS- 1719357, the Simons Fellows Program and by the  Infosys Visiting Chair position at the Indian Institute of Science.}
\email{chari@math.ucr.edu}

\address{Department of Mathematics\\ 
  University of California, Riverside\\ 
  900 University Ave., Riverside, CA 92521}
\email{jdavis@math.ucr.edu}

\address{Department of Mathematics\\ Ithaca College\\ 953 Danby Rd., Ithaca, NY, 14850}
\email{rmoruzzi@ithaca.edu}

\dedicatory{To Kolya Reshetikhin, on his 60th birthday.}

\begin{abstract}
The goal of this paper is to understand the  graded  limit of a  family of irreducible prime   representations of the quantum affine algebra associated to a simply-laced simple Lie algebra $\lie g$.  This   family was introduced in \cite{HL1, HL2} in  the context of  monoidal categorification of cluster algebras.
The graded limit of a member of this family is an  indecomposable graded module for the current algebra $\lie g[t]$; or equivalently a module for  the maximal standard parabolic subalgebra in  the  affine Lie algebra $\widehat{\lie g}$. 
In the case when $\lie g$ is of type $A_n$ the problem was studied in \cite{BCM}, where it was shown that   the graded limit is  isomorphic to a level two Demazure module. In this paper we study the case when $\lie g$ is of type $D_n$.  We show that in certain cases the limit is a generalized Demazure module, i.e., it  is  a submodule of a tensor product of level one Demazure modules. We give a presentation of these modules and   compute their  graded character (and hence also the character of the prime representations)  in terms of Demazure modules of level two.
\end{abstract}
  
    \maketitle

\section*{Introduction}
 Kirillov--Reshetikhin modules (or KR-modules) are a family of finite-dimensional modules for the  quantum affine algebra associated to a simple Lie algebra.  These modules were introduced and originally studied by Kirillov and Reshetikhin in \cite{KR} in connection with the closely related Yangians. The interest in these representations arose   from their study of solvable lattice models which allowed them to  formulate an important conjecture on the character of the tensor products of these modules.
Subsequently, with the introduction of $q$-characters in \cite{FR}, which was further developed in \cite{FM1, FM2}, it became more natural and more tractable to study the Kirillov-Reshetikhin modules for the quantum affine algebra associated to a simple Lie algebra. In \cite{Nak1, Nak2} the Kirillov--Reshetikhin conjecture was solved by geometric methods when $\lie g$ is simply-laced. The non-simply laced case was done in \cite{Hkir} using a $q$-character approach. A combinatorial version of this conjecture was proved in \cite{deFK}.
 A crystal basis approach to these modules was first developed in \cite{HKOTT, HKOTY} and  there is a vast literature on crystal bases for these modules (see for instance  \cite{LNSSS, OS}). We emphasize here that our references to the  topics discussed in the introduction are nowhere  near comprehensive, that would require a separate volume!
 
One approach to understanding the character of KR-modules
and more generally the character of  finite-dimensional representations of quantum affine algebras is to study  the classical ($q\to 1$)  limit of these modules.  The systematic study of this approach was begun in  \cite{CPweyl}, where  a necessary and sufficient condition for the  existence of the limit was proved.  All  the interesting representations studied in the literature  admit a classical limit. The limit, when it exists,  is  a  finite-dimensional module for the corresponding affine Lie algebra and so also a   module for current algebra $\lie g[t]$ which is the Lie algebra  of polynomial maps from $\mathbb C\to \lie g$.
 The notion of the graded limit  of a representation of a quantum affine algebra was
developed in \cite{Cha01,CMkir}. We refer the reader to \cite{Kedem} for a  discussion  of the role of graded limits in the KR-conjecture.  	  
An interesting problem is to give a presentation of the graded limit and to compute its graded character; this is particularly so, because a presentation and the character of an irreducible finite-dimensional module for the quantum affine algebra is not known in general. This problem was studied for a class of representations called minimal affinizations in \cite{moura, Naoi2, Naoi3} when $\lie g$ is of classical type and in \cite{lin,moura2} for some exceptional types.

 D. Hernandez and B. Leclerc introduced in \cite{HL1, HL2} the notion of a monoidal categorification of cluster algebras.
 They defined certain  tensor subcategories of representations of the quantum affine algebra and showed that their Grothendieck ring was a cluster algebra of Dynkin type. The cluster variables are particular
  \lq prime\rq\ irreducible modules (which we shall refer to as HL-modules) of the quantum affine algebra  and a cluster monomial is an irreducible tensor product of HL-modules.  (Recall that a representation is said to be prime if it is not isomorphic to a tensor product of non-trivial representations).  The HL-modules   include some of the well-known ones, such as certain KR-modules,  but also include  many new examples of non-isomorphic prime irreducible representations which are not well-understood.
  
 We are interested in the graded limit of HL-modules. In the case of $A_n$ this problem was studied in \cite{BC, BCM} and we discuss this briefly. The results of  those papers imply  that the graded limit of an HL-module is isomorphic to a Demazure module occurring in  a level two representation of the affine Lie algebra.
  The graded character of Demazure modules is known to be given by the Demazure character formula. In recent work \cite{BCSW},   the graded characters of HL-modules in type $A_n$ are  given explicitly in terms of  Macdonald polynomials. The appearance of level one Demazure modules as graded limits goes back to the work of \cite{CL, FL, Naoi} and we say no more about it here.

In this paper we are interested in the limits of HL--modules in type $D_n$. We restrict our attention to the modules which correspond to the   roots of $D_n$ in which the simple roots occur with multiplicity at most one.
Even with this restriction, the situation is substantially more complicated than  $A_n$. We prove that the graded limit of an HL-module  is a  generalized Demazure module, i.e. it occurs  as a particular submodule of the tensor product of level one Demazure modules. The study of generalized Demazure modules is not very well-developed and in particular, no presentation or character formula is known in general for these modules. 
Thus,  in the first two sections we work entirely in the context of the current algebra $\lie g[t]$ where $\lie g$ is of type $D_n$ and give a presentation and a graded character formula for  a family of generalized Demazure modules. In the final section of the paper we explain briefly the connection with the work of Hernandez and Leclerc.  These connections are similar to the $A_n$ case which is explained in great detail in \cite{BC, BCM}.  There are  difficulties which  arise when we take the HL-modules  associated with the remaining  cluster variables; namely those  which involve  a  simple root with multiplicity two. However, it looks likely that these will have a flag where the successive quotients are level two Demazure modules and we hope to study this elsewhere.

\section{Generalized Demazure Modules}\label{s1}
In this section we  set up the basic notation and define  the Demazure and generalized Demazure modules associated with an  integrable highest weight representation of the untwisted affine Lie algebra  $D_n^{(1)}$.  In the latter part of the section we state and prove the main result assuming two key steps in the proof.
    \subsection{}
        Throughout this paper $\lie{g}$ will denote a simple Lie algebra of type $D_n$ and   $\lie h$  a fixed Cartan subalgebra.
        Let $\kappa:\lie g\times\lie g\to\mathbb C$ be the Killing form on $\lie g$ and $(\ ,\ ):\lie h^*\times \lie h^*\to\mathbb C$ be the induced symmetric bilinear form on $\lie h^*$.
        Let $R$ and  $W$ denote  the set of roots and the Weyl group respectively of the pair $(\lie{g},\lie{h})$. Fix a set of simple roots
         $\{\alpha_i : 1\le i\le n\}$   where $\alpha_{n-1}, \alpha_n$ are the spin nodes and $\alpha_{n-2}$ is the trivalent node and 
        fix  a set of fundamental weights $\{\omega_i : 1\le i\le n\}$ satisfying $(\omega_i,\alpha_j)=\delta_{i,j}$. It will be convenient to set $\omega_0=0=\omega_{n+1}$.
        As usual we let $Q$ and $P$ be the integer span of the simple roots and fundamental weights respectively; the sets $Q^+$ and $P^+$ are defined in the obvious way and we set $R^+=R\cap Q^+$. 
     
   Set,
  \begin{gather*}\alpha_{i,j}=\alpha_i+\cdots+\alpha_j, \ \ \ \ 
  \alpha_{i,n}=\alpha_i+\cdots+\alpha_{n-2}+\alpha_n,  \ \ 1\le i\le j\le n-1,\\
   \beta_{i,j}=\alpha_i+\cdots+\alpha_{j-1}+2(\alpha_j+\cdots+\alpha_{n-2})+\alpha_{n-1}+\alpha_n,\ \ 1\le i<j\le n-1,\\
   h_{i,j}=h_i+\cdots +h_j,\ \ 1\le i\le j\le n-1,\\ h_{i,n}=h_i+h_{i+1}+\cdots+h_{n-2}+h_n,\ \ 1\le i\le n-2.\end{gather*} 
   We shall  adopt the convention that $\alpha_{i,j}=h_{i,j}=0$ if $i>j$ or if $(i,j)=(n-1,n)$.  Note that $$R^+=\{\alpha_{i,j}: 1\le i\le j\le n\}\sqcup\{\beta_{i,j}: 1\le i<j\le n-1\}.$$
   Define a partial order on $P$ by $\lambda \le \mu$ iff $\mu - \lambda \in Q^+$. 
      Finally, let $\{x^\pm_\alpha, h_i: \alpha\in R^+, i\in I\}$ be a Chevalley basis of $\lie g$ and set  $x_{\alpha_i}^\pm =x_i^\pm $. Let $\lie g=\lie n^-\oplus\lie h\oplus\lie n^+$ be the corresponding triangular decomposition.
      \subsection{}
      
      Let $\widehat{\lie g}$ be the untwisted affine Lie algebra associated to $\lie g$ and recall that it has the following explicit realization.
      Given an indeterminate $t$,  let $\Bbb{C}[t,t^{-1}]$ be the  algebra of Laurent polynomials and set
$$\widehat{\lie{g}} = \lie{g} \otimes \Bbb{C}[t,t^{-1}]\oplus \Bbb{C}c \oplus \Bbb{C} d.$$ Define the commutator on $\widehat{\lie g}$ by requiring $c$ to be central, 
       and $$[x\tensor t^r , y\tensor t^s ] = [x,y]\tensor t^{r+s} + r\delta_{r+s,0}\kappa(x,y)c, \ \ [d,x\tensor t^r] = r(x\tensor t^r), \ \ x , y \in \lie{g}, \ \ r,s\in \Bbb{Z}.$$
       Then $$\widehat{\lie h}=\lie h\oplus\mathbb Cc\oplus\mathbb Cd $$ is a Cartan subalgebra of $\widehat{\lie g}$ and we  extend an element $\lambda\in \lie h^\ast$ to an element of $\widehat{\lie h}^\ast$ by setting $\lambda(c)=0=\lambda(d)$. Define  elements $\Lambda_0, \delta,\alpha_0\in\widehat{\lie h}^\ast$ by $$\Lambda_0(\lie h\oplus\mathbb Cd)=0,\ \ \Lambda_0(c)=1,\ \ \delta(\lie h\oplus\mathbb Cc)=0,\ \  \delta(d)=1,\ \ \alpha_0=-{\beta_{1,2}}+\delta.$$  
       The set $\{\alpha_i: 0\le i\le n\}$ is a     set of simple roots for the pair $(\widehat{\lie g}, \widehat{\lie h})$ and the corresponding set of affine  fundamental  weights is given by  $\{\Lambda_i: 0\le i\le n\}$ where $$\Lambda_i=\omega_i+a_i^\vee\Lambda_0,\ \  \  a_i^\vee =2, \ \ 2\le i\le n-2, \ \ a_i^\vee =1,\ \ i=1,n-1,n.$$
       Let $\widehat W$ be the affine Weyl group and note that it  contains an isomorphic copy of $W$.
      Define the affine root lattice $\widehat Q$ and affine weight lattice $\widehat P$ and the subsets $\widehat Q^+$, $\widehat P^+$ and the set of affine positive roots in the natural way. 
       Let $$\widehat {\lie b}=\widehat{\lie h}\oplus (\lie n^+\otimes 1)\oplus (\lie g\otimes  t\mathbb C[t]),\ \ \widehat{\lie p}= (\lie n^-\otimes 1)\oplus\widehat{\lie b} $$
       be the Borel subalgebra and the standard  maximal parabolic subalgebra of $\widehat{\lie g}$ with respect to $\widehat R^+$. 
      Let $\widehat{\lie n}^+$ be the  commutator subalgebra of $\widehat{\lie b}$ and note that the commutator subalgebra of $\widehat{\lie p}$ is just the current algebra  $\lie g[t]= \lie g\otimes \mathbb{C}[t]$ associated to  $\lie g$. The action of $d$ induces
       a canonical  $\Bbb{Z}_+$- grading  on $\lie g[t]$ and on  its universal enveloping algebra with the grade of an element $x\otimes t^r\in\lie g[t]$ being $r$.
       
       \subsection{} Given $\Lambda\in\widehat P^+$ let $\widehat V(\Lambda)$ be the integrable irreducible  $\widehat{\lie g}$-module generated by an element $\widehat v_\Lambda$ with relations $$\widehat{\lie n}^+ \widehat v_\Lambda= 0, \ \ h\widehat v_\Lambda=\Lambda(h)\widehat v_\Lambda,\ \ h\in\widehat{\lie h},$$
       $$(x_i^-\otimes 1)^{\Lambda(h_i)+1} \widehat v_\Lambda=0,\ \ 1\le i\le n, \ \ \ (x_{\beta_{1,2}}^+\otimes t^{-1})^{\Lambda(h_0)+1}\widehat v_\Lambda=0.$$
       We have $$\widehat V(\Lambda) =\bigoplus_{\Lambda'\in \widehat P} \widehat V(\Lambda)_{\Lambda'},\ \  \widehat V(\Lambda)_{\Lambda'}=\{v\in \widehat V(\Lambda): hv=\Lambda'(h)v,\ \ h\in\widehat{\lie h}\}.$$
       It is known that $$\dim \widehat V(\Lambda)_{w\Lambda}= 1,\ \ {\rm{for\ all}}  \ w\in\widehat W,$$ and we let $\widehat V_w(\Lambda)$ be the  $\widehat{\lie b}$-module generated by $\widehat V(\Lambda)_{w\Lambda}$.
       The modules $\widehat V_w(\Lambda)$ are finite-dimensional  and are called the Demazure modules associated to $\Lambda$. Moreover   $\widehat v_\Lambda\in \widehat V_w(\Lambda)$ for all $w\in \widehat W$.  The following is now essentially immediate.
       \begin{lem}\label{hom1} Suppose that $w,w'\in\widehat W$ and $\Lambda\in\widehat P^+$. Then $$\dim\Hom_{\widehat{\lie b}}(\widehat V_w(\Lambda),\widehat V_{w'}(\Lambda))\le 1,$$ and any non-zero  element of this space 
       is injective.
       \hfill\qedsymbol\end{lem}
       
       \subsection{} Given $\Lambda,\Lambda'\in \widehat P^+$ and $w,w'\in\widehat W$ we set $$\widehat V_{w,w'}(\Lambda,\Lambda')=\bold  U(\widehat{\lie b})(\widehat V(\Lambda)_{w\Lambda}\otimes \widehat V(\Lambda')_{w'\Lambda'}).$$ These (along with an obvious extension to an arbitrary number of dominant integral weights and Weyl group elements) are  called the generalized Demazure modules. Further, since $$\dim\Hom_{\widehat{\lie g}}(\widehat V(\Lambda)\otimes \widehat V(\Lambda'),\  \widehat V(\Lambda+\Lambda'))=1$$ it is easy to see that   $$\widehat V(\Lambda)_{w\Lambda}\otimes \widehat V(\Lambda')_{w\Lambda'}\cong \widehat V(\Lambda+\Lambda')_{w(\Lambda+\Lambda')}\ \ {\rm{and\ so}}\ \ 
       \widehat V_{w,w}(\Lambda,\Lambda')\cong \widehat V_w(\Lambda+\Lambda').$$
In this paper we shall be interested in certain special families of  Demazure and generalized Demazure modules; namely those associated to  pairs $(w,\Lambda)$
      such that the restriction of $w\Lambda$ to $\lie h$ is in $(-P^+)$. This restriction guarantees that $\widehat V_w(\Lambda)$ is stable under $\widehat{\lie p}^+$. Equivalently, we can (and do) regard these modules as $\mathbb Z$-graded modules for $\lie g[t]$ and  call them the stable Demazure modules.
      We recall a result proved independently in  \cite{K, M}.
\begin{thm}\label{kumar} Let $\Lambda, \Lambda' \in \widehat{P}^+$ and $w,w'\in\widehat W$. If  $\sigma\in\widehat W$ is such that $\sigma(w\Lambda+w'\Lambda')\in \widehat P^+$, then there exists a projection of $\widehat{\lie g}$-modules $$\widehat V(\Lambda)\otimes \widehat V(\Lambda')\longrightarrow \widehat V(\sigma(w\Lambda+w'\Lambda'))\to 0,$$ which induces an isomorphism of the one-dimensional spaces $$\widehat V(\Lambda)_{w\Lambda}\otimes \widehat V(\Lambda')_{w'\Lambda'}\cong\widehat V(\sigma(w\Lambda+w'\Lambda'))_{w\Lambda+w'\Lambda'}.$$
\hfill\qedsymbol\end{thm}
     The following corollary is immediate from  the preceding theorem and is needed later in this section.  
 \begin{cor}  Retain the notation of the theorem. Assume that $\widehat V_w(\Lambda)$ and $\widehat V_{w'}(\Lambda')$ are stable Demazure modules. Then $\widehat V_{\sigma^{-1}}(\sigma(w\Lambda+w'\Lambda'))$ is also a stable Demazure module and
 there exists a surjective map of $\widehat{\lie p}$-modules  $$\widehat V_{w,w'}(\Lambda,\Lambda')\to\widehat V_{\sigma^{-1}}(\sigma(w\Lambda+w'\Lambda'))\to 0.$$ \hfill\qedsymbol
 \end{cor}

 \subsection{} Let $w_\circ$ be the longest element in $W$. Suppose that $(w,\Lambda)$ and $(w',\Lambda')$ are such that $\Lambda(c)=\Lambda'(c)=\ell$ and $\lambda=w_\circ w\Lambda|_{\lie h}=w_\circ w'\Lambda'|_{\lie h}\in P^+$.   Then it is known \cite{FL, Naoi} that $V_w(\Lambda)$ is isomorphic to $V_{w'}(\Lambda')$ as $\lie g[t]$-modules up to an overall shift of the action of $d$.  In particular we can regard the stable Demazure modules as being indexed by a pair $(\ell,\lambda)\in\mathbb N\times P^+$.
 
 In what follows we prefer to regard the stable (generalized) Demazure modules as $\mathbb Z$-graded modules for the current algebra of $\lie g[t]$, where the grading is given by the action of $d$. The maps between graded modules will be of degree zero. Given $s\in\mathbb Z$ and a graded $\lie g[t]$-module we denote by $\tau_s^*V$ the graded  $\lie g[t]$--module obtained by shifting the grades  by $s$.

      We recall a  presentation of the stable Demazure modules  given in \cite{CV}.
      Let $D(\ell,\lambda)$ be the $\lie g[t]$-module generated by an element $w_\lambda$ with the following defining relations: for $1\le i\le n$, $\alpha\in R^+$, $h\in\lie h$, we have
     \begin{equation}\label{localweyl} \lie n^+[t] w_\lambda=0,\ \ (h\otimes t^r)w_\lambda=\delta_{r,0}\lambda(h)w_\lambda,\ \ (x_i^-\otimes 1)^{\lambda(h_i)+1}w_\lambda=0,\end{equation}
     \begin{equation}\label{deml} (x_\alpha^-\tensor t^{s_\alpha})w_\lambda =0,\ \  (x_\alpha^-\tensor t^{s_\alpha - 1})^{m_\alpha + 1}w_\lambda=0,   
      \end{equation} where we write $\lambda(h_\alpha)=\ell (s_\alpha-1)+m_\alpha$ with $s_\alpha\in\mathbb N$ and $0<m_\alpha\le \ell$. The $\mathbb Z$-grading on $D(\ell,\lambda)$ is given by requiring $w_\lambda$ to have grade zero. 
      
      In the case  $\ell=1,2$  the relations can be further streamlined. The following was proved when $\ell=1$ in \cite{CPweyl} and when $\ell=2$ in \cite{BCM}.
      \begin{lem}\label{simplfiy} If $\ell=1$ the relations in \eqref{deml} are a consequence of the relations in \eqref{localweyl}.  If $\ell=2$ then the second relation in \eqref{deml} is a consequence of the others.
      \hfill\qedsymbol
      
      \end{lem}

   \subsection{}      In the  language of $\lie g[t]$-modules developed above,  a generalized Demazure module corresponds to considering a tensor product $\tau^*_s D(\ell,\lambda)\otimes \tau_{s'}^*D(\ell',\lambda')$ and taking the $\lie g[t]$-module through $w_\lambda\otimes w_{\lambda'}$. 
      Presentations of generalized Demazure modules have not been much studied although there is some work on the subject.  They  play an important role in the study of classical limits of representations of quantum affine algebras and first arose in this context in \cite{Naoi2, Naoi3} (see \cite{ravinder} for other examples).
     In this paper we shall give a presentation of a particular family of generalized Demazure modules, namely those of the form $$D(\lambda,\mu):=\bold U(\lie g[t])(w_\lambda\otimes w_\mu)\subset D(1,\lambda)\otimes D(1,\mu), $$  with some restrictions on the pair $(\lambda,\mu)\in P^+\times P^+$.  Our motivations for considering this problem was discussed extensively in the introduction and further details can also be found in Section \ref{s4}.  We note the  following reformulation of Corollary \ref{kumar}.
     \begin{lem}\label{usekumar} There exists a (unique up to scalars)  map $\eta_{\lambda,\mu}: D(\lambda,\mu)\to D(2,\lambda+\mu)\to 0,$ of  $\lie g[t]$-modules  extending the assignment $w_\lambda\otimes w_\mu\to w_{\lambda+\mu}$.\hfill\qedsymbol
     \end{lem}

      \subsection{}  Given $\lambda,\mu\in P^+$, let $V(\lambda,\mu)$ be the $\lie g[t]$-module generated by an element $w_{\lambda,\mu}$ satisfying the following defining relations: 
      \begin{gather}\label{localweyla} \lie n^+[t] w_{\lambda,\mu}=0,\ \ (h\otimes t^r)w_{\lambda,\mu}=\delta_{r,0}(\lambda+\mu)(h)w_{\lambda, \mu},\ \ (x_i^-\otimes 1)^{(\lambda+\mu)(h_i)+1}w_{\lambda,\mu}=0,\\ \label{gdem} (x_\alpha^-\tensor t^{\max\{\lambda(h_\alpha),\mu(h_\alpha)\}})w_{\lambda,\mu} =0,\end{gather} where $1\le i\le n$, $h\in\lie h$ and $\alpha\in R^+$.
      Define a  grading on $V(\lambda,\mu)$ by requiring the grade of $w_{\lambda,\mu}$ to be zero. Clearly $V(\lambda,\mu)$ is a graded quotient of $D(1,\lambda+\mu)$ and $V(\lambda,\mu)\cong V(\mu,\lambda)$. Moreover Lemma \ref{simplfiy} shows  that if $\mu=0$ then $V(\lambda,0) \cong D(1,\lambda).$
     The following  is immediate from equations  \eqref{localweyl}-\eqref{gdem} and Lemma \ref{usekumar}.
     \begin{lem}\label{lem:VtoD}
     The assignments $w_{\lambda,\mu}\to w_{\lambda+\mu}$ and $w_{\lambda,\mu}\to w_\lambda\otimes w_\mu$ define   surjective maps of graded  $\lie g[t]$-modules $$\psi_{\lambda,\mu}: V(\lambda,\mu)\to D(2,\lambda+\mu),\ \ \varphi_{\lambda,\mu}: V(\lambda,\mu)\to D(\lambda,\mu)\to 0,  $$ 
     and $\psi_{\lambda,\mu}=\eta_{\lambda,\mu}\circ\varphi_{\lambda,\mu}.$
     \hfill\qedsymbol
     \end{lem}
     \subsection{} 
     Let
     $$P^+(1) = \{ \lambda \in P^+ : \lambda(h_i) \le 1, \ \ 1 \le i \le n\}.$$
     It is obvious that any  $\lambda\in P^+(1)$ can be written  uniquely (up to order) as a sum $\lambda=\lambda_1+\lambda_2$ where $\lambda_s\in P^+(1)$, $s=1,2$ and satisfy the following:  for $1\le i \le j\le n$ and $(i,j)\ne (n-1,n)$ 
   $$\lambda_r(h_i)=1=\lambda_r(h_j)\ \implies \lambda_p(h_s)=1\ {\rm{for\ some}}\ \ i<s<j,\ \ \{r,p\}=\{1,2\},$$  and $$\lambda(h_{n-1}+h_n)>0\implies \lambda_p(h_{n-1}+h_n)=0 \ \ {\rm{for\ some}}\ \ p\in\{1,2\}.$$
   It will be convenient to call    $(\lambda_1, \lambda_2)\in P^+\times P^+$ an {\em interlacing} pair if  $\lambda_1+\lambda_2\in P^+(1)$ and the preceding two conditions hold.
  
  \noindent {\bf Examples.}
   The pairs $(\omega_i,0)$ for $0\le i\le n$, $(\omega_{n-1}+\omega_{n}, 0)$ and the elements of the set
  $\{(\omega_i,\omega_j): 0\le i\ne j\le n\ \ (i,j)\ne (n-1,n)\}$ are interlacing.
 The pair  $(\omega_2+\omega_5,\omega_3+\omega_6+\omega_7)$ is  interlacing in $D_7$ but not in $D_8$.

 \subsection{}  The following is the main result of this note. 
   \begin{thm}\label{main} Let $(\lambda_1,\lambda_2)\in P^+\times P^+$ be an interlacing pair. For all $\nu\in P^+$ the map $$\varphi_{\lambda_1+\nu,\lambda_2+\nu}: V(\lambda_1+\nu,\lambda_2+\nu)\to D(\lambda_1+\nu,\lambda_2+\nu),$$ is an isomorphishm. Moreover  $D(\lambda_1+\nu,\lambda_2+\nu)$ admits a flag and the quotients are isomorphic to  $\tau_s^* D(2,\mu)$ for some $s\in\mathbb Z_+$, $\mu\in P^+$. 
   \end{thm}
    \begin{rem} One can give a precise formula for the graded character of $D(\lambda_1+\nu,\lambda_2+\nu)$ and this can be found in Section \ref{chgr}.
    \end{rem}
  
  \begin{rem}  An analogous result  was proved for $A_n$ in \cite{BCM} when $\nu=0$.  In fact in that case it was proved that $$V(\lambda_1,\lambda_2)\cong D(\lambda_1,\lambda_2)\cong D(2,\lambda_1+\lambda_2).$$ Our methods are different and work for $A_n$ as well and prove the  more general case when $\nu\ne 0$.
   \end{rem}
   
  \begin{rem} 
   It is plausible that the map $\varphi_{\lambda,\mu}: V(\lambda,\mu)\to D(\lambda,\mu)$ is always an isomorphism. This is known to be true for  $\lie{sl}_2$  and examples exist in $\lie{sl}_n$ (see for instance \cite{BP, Naoi4}).
   
   \end{rem}
   {\em In the rest of the section we shall assume that $(\lambda_1,\lambda_2)$ is an interlacing pair and  $\lambda=\lambda_1+\lambda_2$. Moreover the property of interlacing pairs means that we can and  will assume without loss of generality that if $\lambda(h_{n-1}+h_n)>0$ then $\lambda_1(h_{n-1})>0$. This means that if 
   $p$ is maximal such that $\lambda(h_{p+1})=1$ then $p\le n-2$ and $\lambda_1(h_{p+1})=1$ and  $\lambda_2(h_{p+1})=0$.}
   \subsection{} The theorem is proved in several steps. The first reduction is the following proposition which gives a condition for the generalized Demazure module to be isomorphic to a Demazure module.
   \begin{prop} \label{demiso} If $\lambda(h_{n-1}+h_n)= 1$ or if  $\lambda=\omega_{i-1}+\omega_i+\delta_{i+1,n}\omega_n$ for $1\le i\le n-1$ then for all $\nu\in P^+$ we have$$V(\lambda_1+\nu,\lambda_2+\nu)\cong D(2, 2\nu+\lambda)\cong D(\lambda_1+\nu,\lambda_2+\nu).$$
   \end{prop}
   \subsection{}\label{defbeta}  Suppose that  $\lambda(h_{n-1}+h_n)\in\{0,2\} $ and   $\lambda\ne\omega_{i-1}+\omega_i+\delta_{i+1,n}\omega_n$ and let $p\le n-2$ be maximal such that $\lambda_1(h_{p+1})=1=\lambda(h_{p+1})$. Let  $1\le p'\le p\le n-2$ be maximal so that $(\lambda_1-\lambda_2)(h_{p',p})=0$. Observe that if $p'<p$ then the definition of interlacing pairs forces $\lambda_2(h_p)=1=\lambda_1(h_{p'}).$
   Set $$\beta_\lambda=\beta_{p',p+1}=\alpha_{p'}+\cdots+\alpha_p+2(\alpha_{p+1}+\cdots+\alpha_{n-2})+\alpha_{n-1}+\alpha_n$$ and note that  $\lambda_1(h_{\beta_\lambda})=3-\delta_{p',p}$ and $ \lambda_2(h_{\beta_\lambda})= (1-\delta_{p'p}).$ It is helpful to note that the assumption $\lambda(h_{n-1}+h_n)\in\{0,2\} $ implies that $p=n-2$ iff $\lambda(h_{n-1}+h_n)=2$.
   
    \begin{lem}\label{diff2proved} Suppose that  $\lambda(h_{n-1}+h_n)\in\{0,2\} $ and   $\lambda\ne\omega_{i-1}+\omega_i+\delta_{i+1,n}\omega_n$.
  Then $\lambda_1-\beta_\lambda\in P^+$ and there exists $\nu_0\in P^+$ such that     $  (\lambda_1 - \beta_\lambda - \nu_0,\lambda_2-\nu_0)$ is an interlacing pair. \end{lem}

   \begin{pf}  It is simple to see  that with our assumptions $$\lambda_1-\beta_\lambda=\lambda_1-\omega_{p+1}-\delta_{p,n-2}\omega_n- (1-\delta_{p',p})\omega_{p'}+\omega_{p'-1}+(1-\delta_{p',p})\omega_p\in P^+. $$   
   Taking $\nu_0 = \lambda_2(h_{p'-1})\omega_{p'-1} +  (1-\delta_{p',p})\lambda_2(h_p)\omega_p$  it is easy to check that  $ (\lambda_1 - \beta_\lambda - \nu_0,\lambda_2-\nu_0)$ is an  interlacing pair.

   \end{pf}
   
   \subsection{} The second reduction is the following proposition which in particular,  gives an upper bound for the dimension of $V(\lambda_1+\nu,\lambda_2+\nu)$.
   \begin{prop}\label{res} Suppose that  $\lambda(h_{n-1}+h_n)\in\{0,2\} $ and   $\lambda\ne\omega_{i-1}+\omega_i+\delta_{i+1,n}\omega_n$. There exists a right exact sequence of $\lie g[t]$-modules \begin{equation} \tau^\ast_{(\lambda_1+\nu)(h_{\beta_\lambda})-1} V(\lambda_1+\nu-{\beta_\lambda},\lambda_2+\nu)\rightarrow V(\lambda_1+\nu,\lambda_2+\nu)\xrightarrow{} D(2,\lambda+2\nu)\to
   0,\end{equation} with $w_{\lambda_1+\nu-{\beta_\lambda},\lambda_2+\nu}\to (x_{\beta_\lambda}^-\otimes t^{(\lambda_1+\nu)(h_{\beta_\lambda})-1})w_{\lambda_1+\nu,\lambda_2+\nu}$.
    \end{prop}

   \subsection{} Assuming Proposition \ref{demiso} and Proposition \ref{res} we complete the proof of Theorem \ref{main}.
   The proof is by an induction  with respect to  the partial order on $P^+$. The minimal elements with respect to this order are $0, \omega_1,\omega_{n-1},\omega_n$ and 
   Proposition \ref{demiso} shows that induction begins.  Moreover  the  proposition also establishes the theorem   if $\lambda=\omega_{i-1}+\omega_{i}+\delta_{i+1,n}\omega_n$ or $\lambda(h_{n-1}+h_n)=1$. Hence we have only to prove   the inductive step when $\lambda(h_{n-1}+h_n)\in\{0,2\}$ and $\lambda\ne\omega_{i-1}+\omega_{i}+\delta_{i+1,n}\omega_n$.
   We shall need the following result to complete the proof of the inductive step.
   \begin{lem}\label{weylinc} 
 If  $\lambda(h_{n-1}+h_n)\in\{0,2\}$ and $\lambda\ne\omega_{i-1}+\omega_{i}+\delta_{i+1,n}\omega_n$, there exists an inclusion of $\lie g[t]$-modules \begin{equation} \tau^\ast_{(\lambda_1+\nu)(h_{\beta_\lambda})-1}D(1,\lambda_1+\nu-{\beta_\lambda})\hookrightarrow D(1,\lambda_1+\nu),\end{equation} which sends to $w_{\lambda_1+\nu-{\beta_\lambda}}\to(x_{\beta_\lambda}^-\otimes t^{(\lambda_1+\nu)(h_{\beta_\lambda})-1})w_{\lambda_1+\nu}. $  
    \end{lem}
   \begin{pf} By Lemma \ref{simplfiy} it suffices to prove that $w:=(x_{\beta_\lambda}^-\otimes t^{(\lambda_1+\nu)(h_{\beta_\lambda})-1})w_{\lambda_1+\nu}$ satisfies the relations in \eqref{localweyl}. If $\beta_\lambda-\alpha_i\notin R^+$ then it is clear that $(x_i^+\otimes t^r)w=0$ for all $r\in\mathbb Z_+$. Otherwise we must have $i=p+1$, or $i=p'$ if $p'<p$ or  $i=n$
   if $p=n-2$ and we have to prove that
   $$(x^-_{\beta_\lambda-\alpha_i}\otimes t^{(\lambda_1+\nu)(h_{\beta_\lambda})-1+r})w_{\lambda_1+\nu}=0, \ \ r\in\mathbb Z_+.$$  But this follows from \eqref{deml} since  in all cases $\lambda_1(h_i)=1$ and so $(\lambda_1+\nu)(h_{\beta_\lambda})-1\ge(\lambda_1+\nu)(h_{\beta_\lambda-\alpha_i})$. The second relation in \eqref{localweyl} is trivial if $r=0$ and if $r>0$ then it follows from \eqref{deml} since $$(h\otimes t^r)w=(x_{\beta_\lambda}^-\otimes t^{(\lambda_1+\nu)(h_{\beta_\lambda})+r-1})w_{\lambda_1+\nu}=0.$$ The third relation is immediate since the modules are all finite-dimensional.
   
   \end{pf}
   
   \subsection{}   Lemma \ref{usekumar}, Proposition \ref{res}  and the inductive hypothesis establish the following inequalities:
   \begin{gather*}
       \dim D(\lambda_1+\nu,\lambda_2+\nu)\le \dim V(\lambda_1+\nu,\lambda_2+\nu),\\ \dim V(\lambda_1+\nu,\lambda_2+\nu)\le \dim D(2,2\nu+\lambda)+ \dim D(\lambda_1+\nu-{\beta_\lambda},\lambda_2+\nu).
   \end{gather*}
   The inductive step follows if we prove that
   \begin{equation}\label{dimd} \dim D(\lambda_1+\nu,\lambda_2+\nu)=\dim D(2,2\nu+\lambda)+ \dim D(\lambda_1+\nu-{\beta_\lambda},\lambda_2+\nu).\end{equation} 
   For this,  we observe that Lemma \ref{weylinc} gives an inclusion$$0\to D(1,\lambda_1+\nu-{\beta_\lambda})\otimes D(1,\lambda_2+\nu)\to D(1,\lambda_1+\nu)\otimes D(\lambda_2+\nu), $$ which sends $$w_{\lambda_1+\nu-{\beta_\lambda}}\otimes w_{\lambda_2+\nu}\to ((x_{\beta_\lambda}^-\otimes t^{(\lambda_1+\nu)(h_{\beta_\lambda})-1)})w_{\lambda_1+\nu})\otimes w_{\lambda_2+\nu}.$$ Since $(\lambda_1+\nu)(h_{\beta_\lambda})-1\ge (\lambda_2+\nu)(h_{\beta_\lambda})$  we see that  \eqref{deml} gives $$(x_{\beta_\lambda}^-\otimes t^{(\lambda_1+\nu)(h_{\beta_\lambda})-1)})\left (w_{\lambda_1+\nu}\otimes w_{\lambda_2+\nu}\right)=  \left((x_{\beta_\lambda}^-\otimes t^{(\lambda_1+\nu)(h_{\beta_\lambda})-1})w_{\lambda_1+\nu}\right)\otimes w_{\lambda_2+\nu}.$$
\noindent It follows that we have an inclusion
$$ \iota: D(\lambda_1+\nu-{\beta_\lambda},\lambda_2+\nu)\hookrightarrow D(\lambda_1+\nu,\lambda_2+\nu)$$
and it suffices to prove that the corresponding quotient is isomorphic to $D(2, 2\nu+\lambda)$. By Lemma \ref{lem:VtoD} we have the following surjective maps $$V(\lambda_1+\nu,\lambda_2+\nu)\twoheadrightarrow D(\lambda_1+\nu,\lambda_2+\nu)\twoheadrightarrow  D(2,2\nu+\lambda).$$ These  maps are all unique up to scalars and   Proposition \ref{res} shows that the kernel of the composite map   is generated by the element $(x_{\beta_\lambda}^-\otimes t^{(\lambda_1+\nu)(h_{\beta_\lambda})-1)})w_{\lambda_1+\nu,\lambda_2+\nu}$. Hence  the kernel of  $$D(\lambda_1+\nu,\lambda_2+\nu)\twoheadrightarrow D(2,2\nu+\lambda)$$ is generated by $(x_{\beta_\lambda}^-\otimes t^{(\lambda_1+\nu)(h_{\beta_\lambda})-1)})(w_{\lambda_1+\nu}\otimes w_{\lambda_2+\nu})$. But this means that the latter kernel is precisely the image of $\iota$ and hence the corresponding  quotient is isomorphic to  $ D(2,2\nu+\lambda)$ as needed.

 \subsection{}\label{chgr} Given a graded $\lie g[t]$-module $V=\oplus_{s\in\mathbb Z} V[s]$ we have, $$(\lie g\otimes t^r)V[s]\subset V[r+s],\ \ V_\mu=\bigoplus_{s\in\mathbb Z} V_\mu[s],\ \  V_\mu[s]= V_\mu\cap V[s],$$
 for all $r,s\in\mathbb Z$ and $\mu\in P$. Let $\mathbb Z[v,v^{-1}][P]$ be the group ring of $P$ with coefficients in $\mathbb Z[v,v^{-1}]$ where $v$ is an indeterminate. Let $e(\mu)$, $\mu\in P$ be a basis for the group ring and set $$\ch_{\gr} V=\sum_{\mu\in P}\left(\sum_{ s\in\mathbb Z}\dim V_\mu[s]v^s\right)e(\mu).$$
 Notice that in the course of proving Theorem \ref{main} we also established that  the right exact sequence in Proposition \ref{res} is exact. Hence we get an equality of graded characters \begin{equation}\label{chgr} \ch_{\gr} D(\lambda_1+\nu,\lambda_2+\nu)=\ch_{\gr} D(2,2\nu+\lambda) + v^{(\lambda_1+\nu)(h_{\beta_\lambda})-1}\ch_{\gr} D(\lambda_1-\beta_\lambda+\nu,\lambda_2+\nu).\end{equation}
 Using this we can derive a formula for $\ch_{\gr}D(\lambda_1+\nu,\lambda_2+\nu)$ in terms of Demazure modules as follows. 
 
 Given $\mu\in P^+$ define elements $\lfloor\mu/2\rfloor$ and $\res_2\mu$ of $P^+$ by,
 $$\lfloor\mu/2\rfloor(h_i)=\lfloor\mu(h_i)/2\rfloor, \ \ 1\le i\le n,\ \ \res_2\mu=\mu-2\lfloor\mu/2\rfloor.$$
 Then $\res_2\mu\in P^+(1)$ and we set $\beta_\mu=\beta_{\res_2\mu}$ if $\res_2\mu$ satisfies the conditions of Section \ref{defbeta} and otherwise set $\beta_\mu=0$. Let  $(\mu_1,\mu_2)$ be  the interlacing pair corresponding to $\res_2\mu$ and let  $$r_\mu=\lfloor\mu/2\rfloor(h_{\beta_\mu})+\max\{\mu_1(h_{\beta_\mu}),\mu_2(h_{\beta_\mu})\}-1.$$ 
 Given $\mu\in P^+$ set $$\mu^0=\mu,\ \ \mu^1=\mu^0-\beta_{\mu^0},\dots ,\mu^s=\mu^{s-1}-\beta_{\mu^{s-1}}, $$   where $s$ is minimal so that $\beta_{\mu^s}=0$.
Setting $r_k=r_{\mu^k}$ a repeated application of \eqref{chgr} gives
$$\ch_{\gr}D(\lfloor\mu/2\rfloor+\mu_1, \ \lfloor\mu/2\rfloor+\mu_2)=\sum_{k=0}^s v^{r_k}\ch_{\gr} D(2, \mu^k ).$$

 \section{Proof of Proposition \ref{demiso} and Proposition \ref{res}.}\label{s3}
  In this section  we first establish some technical results on the properties of interlacing pairs and then   prove Proposition \ref{demiso} and   Proposition \ref{res}. 
  
  {\em  Throughout this section we shall assume  without further mention that $(\lambda_1,\lambda_2)$ is an interlacing pair and  $\lambda=\lambda_1+\lambda_2$. We shall   also  assume without loss of generality that there exists $p$  maximal with $\lambda(h_{p+1})=1$ with $p\le n-2$ and  $\lambda_1(h_{p+1})=1$.}
 \subsection{}   \begin{lem} \label{typea}
    For all $1\le i\le j\le n$ we have
   \begin{equation}\label{typear}  |(\lambda_1-\lambda_2)(h_{i,j})|\le 1\end{equation}
   and hence $| (\lambda_1-\lambda_2)(h_\alpha)|\le 2\ \ {\rm{for\ all}}\ \alpha\in R^+.$ Further if $\lambda(h_{n-1}+h_n)=1$ then 
   \begin{equation*}\label{typear} | (\lambda_1-\lambda_2)(h_\alpha)|\le 1\ \ {\rm{for\ all}}\ \alpha\in R^+. \end{equation*}
   \end{lem}
   \begin{pf} Equation \eqref{typear} is proved by induction on $j-i$; induction obviously begins at $i=j$ since $\lambda\in P^+(1)$. If $j>i$
   the inductive step is clear if $\lambda(h_i)=0$ or $\lambda(h_j)=0$. If  $\lambda(h_i)=1=\lambda(h_j)$ then either $\lambda(h_s)=0$ for all $i<s< j$ or there exists $i<p< j$ minimal such that $\lambda_2(h_p)=1$.
  In the first case we have proved that $(\lambda_1-\lambda_2)(h_{i,j})=0$ and in the second case we have proved that $(\lambda_1-\lambda_2)(h_{i,j})=(\lambda_1-\lambda_2)(h_{p+1,j})$ and the inductive step follows.  Since $\beta_{i,j}=\alpha_{i,n}+\alpha_{j,n-1}$ for $1\le i<j\le n-1$ it is immediate that $$|(\lambda_1-\lambda_2)(\beta_{i,j})|\le 2\ \ {\rm{for\ all}}\ 1\le i<j\le n-1.$$ 
  Moreover,  if $\lambda(h_{n-1}+h_n)=1$ our assumptions imply that $\lambda_1(h_{n-1})=1$ and $\lambda_2(h_{n-1}+h_n)=0$ and so
   using  \eqref{typear} we have  
$$-1\le (\lambda_1 - \lambda_2)(h_{i,n}) = (\lambda_1-\lambda_2)(h_{i,n-2})=(\lambda_1-\lambda_2)(h_{i,n-1})-1\le 0, $$ $$1\ge (\lambda_1 - \lambda_2)(h_{j,n-1}) = (\lambda_1-\lambda_2)(h_{j,n-2})+1\ge 0. $$ 
Hence we have $$-1\le (\lambda_1-\lambda_2)(h_{\beta_{i,j}})= (\lambda_1-\lambda_2)(h_{i,n})+(\lambda_1-\lambda_2)(h_{j,n-1}) \le 1 $$ and the proof of the lemma  is complete.
   \end{pf}

  \subsection{}

   Recall our convention that $\alpha_{i,j}= 0$ if $i>j$.  Set $$ R(\lambda_1,\lambda_2)=\{\beta_{i,j}: (\lambda_1-\lambda_2)(h_{\beta_{i,j}})=\pm 2\}.$$
  Note that Lemma \ref{typea} shows that  $R(\lambda_1,\lambda_2)=\emptyset$ if $\lambda(h_{n-1}+h_n)=1$ or $ \lambda= \omega_{i-1}+\omega_i + \delta_{i,n-1}\omega_n.$ The next result establishes the converse.
 \begin{lem} \label{min} Suppose that $ \lambda(h_{n-1}+h_n)\in\{0,2\}$  and $\lambda\ne \omega_{i-1}+\omega_i + \delta_{i,n-1}\omega_n.$
  Then $\beta_\lambda\in R(\lambda_1,\lambda_2)$ and more generally,
 $$  \beta_{i,j}\in R(\lambda_1,\lambda_2)\iff  \beta_{i,j}=\alpha_{i,p'-1}+\beta_\lambda+\alpha_{j,p}
 \ \ {\text{and}}\ \  
  (\lambda_1-\lambda_2)(h_{i,p'-1})=0=(\lambda_1-\lambda_2)(h_{j,p})=0.$$ 
  \end{lem}
 \begin{pf} 
 Recall from Section \ref{defbeta} that $\beta_\lambda=\beta_{p',p+1}$ where $1\le p'\le p$ is maximal with $(\lambda_1-\lambda_2)(h_{p',p})=0$ and that  $(\lambda_1-\lambda_2)(h_{\beta_\lambda})=2$.
 Hence a calculation shows that   $$\beta_{i,j}\in R(\lambda_1,\lambda_2)\implies i\le p'\ \ {\rm{or}}\ \  p'<j\le p+1,$$ which shows that
    $\beta_{i,j}=\alpha_{i,p'-1}+\beta_\lambda+\alpha_{j,p}.$   
    
    Since $\beta_{i,p+1}\in R^+$ if $i<p'$ we have  $$(\lambda_1-\lambda_2)(h_{\beta_{i,p+1}})=(\lambda_1-\lambda_2)(h_{i,p'-1})+2$$  and
    Lemma \ref{typea} forces  $(\lambda_1-\lambda_2)(h_{i,p'-1})\in\{-1,0\}$. If $j< p+1$ we have $\beta_{p',j}=\beta_{p',p+1}+\alpha_{j,p}$ and hence  a similar argument shows that $(\lambda_1-\lambda_2)(h_{j,p})\in\{-1,0\}$. Together with hypothesis that  $(\lambda_1-\lambda_2)(h_{\beta_{i,j}})=\pm 2$ we get
   $(\lambda_1-\lambda_2)(h_{i,p'-1})=0=(\lambda_1-\lambda_2)(h_{j,p})$ as needed. The converse is obvious.
   \end{pf}
  \begin{rem} Notice that in particular we have  established  (with our conventions on $(\lambda_1,\lambda_2)$) that  $\beta\in R(\lambda_1,\lambda_2)$ iff $(\lambda_1-\lambda_2)(h_\beta)=2$.
  \end{rem}

\subsection{} Assume that $\lambda(h_{n-1}+h_n)=1$ or that $\lambda=\omega_{i-1}+\omega_i+\delta_{i,n-1}\omega_n$. By Lemma \ref{typea} we have $(\lambda_1-\lambda_2)(h_\alpha)\in\{-1,0,1\} $ for all $\alpha\in R^+$, and so $$ \max\{(\nu+\lambda_1)(h_\alpha),(\nu+\lambda_2)(h_\alpha)\}=\nu(h_\alpha)+\lceil(\lambda(h_\alpha)/2\rceil.$$ The defining relations give
 $$V(\lambda_1+\nu,\lambda_2+\nu)\cong D(2, 2\nu+\lambda_1+\lambda_2).$$ 
Since the maps in   Lemma \ref{usekumar} and Lemma \ref{lem:VtoD} are unique up to scalars it now follows  that the map $$V(\lambda_1+\nu,\lambda_2+\nu)\twoheadrightarrow D(\lambda_1+\nu,\lambda_2+\nu)\twoheadrightarrow D(2, 2\nu+\lambda_1+\lambda_2),$$ is an isomorphism and hence all the maps are isomorphisms. Proposition \ref{demiso} is proved.

 \subsection{} We turn to the proof of Proposition \ref{res}.  We have $R(\lambda_1,\lambda_2)\ne\emptyset$ and  our conventions imply  that $\lambda_1(h_{\beta_\lambda})= 3 - \delta_{p',p}$ and $ \lambda_2(h_{\beta_\lambda})=1 - \delta_{p',p}$. Using the remark following   Lemma \ref{min}  we have  $\lambda_1(h_\beta)\ge\lambda_2(h_\beta)$ for all $\beta\in R(\lambda_1,\lambda_2)$. 
Hence for all $\beta\in R(\lambda_1,\lambda_2)$, \begin{equation*}\label{maxceila} (\nu+\lambda_1)(h_\beta)=\max\{(\nu+\lambda_1)(h_\beta),(\nu+\lambda_2)(h_\beta)\}=\nu(h_\beta)+\lceil(\lambda(h_\beta)/2\rceil +1,.\end{equation*}
An inspection of the   defining relations of the modules shows  that   the kernel $K$ of the canonical map  $$\psi_{\lambda_1+\nu,\lambda_2+\nu}: V(\lambda_1+\nu,\lambda_2+\nu)\to D(2,\lambda + \nu) \to 0,\ \ w_{\lambda_1+\nu,\lambda_2+\nu}\to w_{2\nu+\lambda},$$
is generated by the elements  
 $$(x_\beta^-\otimes t^{(\nu+\lambda_1)(h_{\beta})-1})w_{\lambda_1+\nu,\lambda_2+\nu}, \ \ \beta\in R(\lambda_1,\lambda_2).$$
  The proof of Proposition \ref{res}  is  now  shown in two steps. 
The first step is to show that $$K=\bold U(\lie g[t])(x_{\beta_\lambda}^-\otimes t^{(\nu+\lambda_1)(h_{\beta_\lambda})-1)})w_{\lambda_1+\nu,\lambda_2+\nu}.$$ 
For this, let $\beta\in R(\lambda_1,\lambda_2)$ and write $\beta=\alpha_{i,p'-1}+\beta_{p',p+1}+\alpha_{j,p}$ as in Lemma \ref{min} and assume that $i\le p'-1$ or $j\le p$ (otherwise $\beta=\beta_\lambda$ and there is nothing to prove). 
The defining relations
  $$(x_{i,p'-1}^-\otimes t^{(\nu+\lambda_1)(h_{i,p'-1})})w_{\lambda_1+\nu,\lambda_2+\nu} = 0 = (x_{j,p}^-\otimes t^{(\nu+\lambda_1)(h_{j,p})})w_{\lambda_1+\nu,\lambda_2+\nu}, $$
  imply if  $i\le p'-1$ and  $j\le p$ \begin{gather*}(x_{i,p'-1}^-\otimes t^{(\nu+\lambda_1)(h_{i,p'-1})})(x_{j,p}^-\otimes t^{(\nu+\lambda_1)(h_{j,p})})(x^-_{\beta_\lambda}\otimes t^{(\nu+\lambda_1)(h_{\beta_\lambda})-1}
)w_{\lambda_1+\nu, \lambda_2+\nu}=\\ [x_{i,p'-1}^-, [x_{j,p}^-, x_{\beta_\lambda}^-]]\otimes t^{(\nu+\lambda_1)(h_{\beta_{i,j}})-1})w_{\lambda_1+\nu,\lambda_2+\nu}= (x_{\beta_{i,j}}^-\otimes t^{(\nu+\lambda_1)(h_{\beta_{i,j}})-1})w_{\lambda_1+\nu,\lambda_2+\nu}.\end{gather*} An obvious modification works if $i\ge p'$ or $j\ge p+1$   and  the first step is proved. 

The next step is to show  that the  element $(x_{\beta_\lambda}^-\otimes t^{(\nu+\lambda_1)(h_{\beta_\lambda})-1})w_{\lambda_1+\nu,\lambda_2+\nu}$ satisfies the defining relations of the element $w_{\lambda_1+\nu-{\beta_\lambda},\lambda_2+\nu 
}\in  V(\lambda_1+\nu-{\beta_\lambda},\lambda_2+\nu)$. This establishes
the existence of the map $V(\lambda_1+\nu-{\beta_\lambda},\lambda_2+\nu)\rightarrow K\to 0$.

\subsection{} We prove that $(x_{\beta_\lambda}^-\otimes t^{(\nu+\lambda_1)(h_{\beta_\lambda})-1})w_{\lambda_1+\nu,\lambda_2+\nu}$ satisfies the relations in \eqref{localweyla} with $(\lambda_1,\lambda_2)$ replaced by $(\lambda_1-\beta_\lambda,\lambda_2)$. Since $(x_i^+\otimes t^r)w_{\lambda_1+\nu,\lambda_2+\nu}=0$ for $1\le i\le n$ and $r\in\mathbb Z_+$ the first relation in \eqref{localweyla} is immediate if $\beta_\lambda-\alpha_i\notin R^+$.  Otherwise  we must prove that
$$(\dagger)\ \  (x_{\beta_\lambda-\alpha_i}^-\otimes t^{(\nu+\lambda_1)(h_{\beta_\lambda})-1+r})w_{\lambda_1+\nu,\lambda_2+\nu}=0.$$ Since $\beta_\lambda=\beta_{p',p+1}$ it follows that either  $i=p+1$ and if $p'<p$ then we can also have   $i=p'$ and if $p=n-2$ then, in addition we can have  $i=n$. In the first two cases we have $\lambda_1(h_i)=1$ and $\lambda_2(h_i)=0$ (see Section \ref{defbeta}). In particular $$(x_{\beta_\lambda-\alpha_i}^-\otimes t^{\nu(h_{\beta_\lambda-\alpha_i})+ 2-\delta_{p,p'}} )w_{\lambda_1+\nu,\lambda_2+\nu}=0, $$ 
is a defining relation and  $(\dagger)$ follows by applying $(h\otimes t^r)$ since $$(\nu+\lambda_1)(h_{\beta_\lambda})-1+r=\nu(h_{\beta_\lambda})+2- \delta_{p,p'}+r\ge \nu(h_{\beta_\lambda-\alpha_i})+2-  \delta_{p,p'},\ \  r\in\mathbb Z_+.$$  A similar argument establishes the case when $p=n-2$ and $i=n$.

 It is trivial to see that the second relation in \eqref{localweyl}  holds and the third is then immediate since $ V(\nu+\lambda_1,\nu+\lambda_2)$ is finite-dimensional.
\subsection{} We need the following result to prove that $(x_{\beta_\lambda}^-\otimes t^{(\nu+\lambda_1)(h_{\beta_\lambda})-1})w_{\lambda_1+\nu,\lambda_2+\nu}$ satisfies the relations in \eqref{gdem}.
\begin{lem} The relation $(x_\alpha^-\otimes t^{\nu(h_\alpha)+\max\{\lambda_1(h_\alpha), \lambda_2(h_\alpha)\}})w_{\lambda_1+\nu,\lambda_2+\nu}=0$ holds in $V(\lambda_1+\nu,\lambda_2+\nu) $ for all $\alpha\in R^+$ iff it holds for all the  elements of the set $$\{\alpha_{i,j}: 1\le i\le j\le n,\ \ \lambda(h_i)=1=\lambda(h_j),
\ \ \lambda(h_{i,j})=2-\delta_{i,j}\}.$$
\end{lem}
\begin{pf} The forward direction of the Lemma is obvious. For the reverse direction we proceed by induction with respect to the partial order on $R^+$ induced by the partial order $\le $ on $P$. If $\alpha=\alpha_i$ for some $1\le i\le n$ there is nothing to prove. For the inductive step choose $1\le r\le n$ so that $\alpha-\alpha_r\in R^+$. If $\lambda(h_r)=0$   then the result follows from $[x_r^-,x_{\alpha-\alpha_r}^-]=x_\alpha^-$ and $$(x_r^-\otimes t^{\nu(h_r)})w_{\lambda_1+\nu,\lambda_2+\nu}=0= (x_{\alpha-\alpha_r}^-\otimes t^{\nu(h_{\alpha-\alpha_r})+\max\{\lambda_1(h_\alpha),\lambda_2(h_\alpha)\}})w_{\lambda_1+\nu,\lambda_2+\nu},$$ where the second equality is an application  of the inductive hypothesis.
Otherwise,  $\lambda(h_r)> 0$ for all $1\le r\le n$ with $\alpha-\alpha_r\in R^+$, in which case   $\alpha$ is  one of the following: \\
$\bullet$ $\alpha_{i,j}$ with $\lambda(h_i)=1=\lambda(h_j)$,\\  $\bullet$ $\beta_{i,j}$ with $i<j-1$ and $\lambda(h_i)=\lambda(h_j)=1$,\\  
 $\bullet$ $\beta_{j-1,j}$ with $\lambda(h_j)=1$
.\\ In the first two cases choose  $i<s\le j$ minimal with $\lambda(h_s)=1$ and $\lambda(h_{i,s})=2$.
The result follows since Lemma \ref{typea} and the inductive hypothesis give $$(x_{i,s}^-\otimes t^{\nu(h_{i,s})+1})w_{\lambda_1+\nu,\lambda_2+\nu}=0=
(x_{\alpha-\alpha_{i,s}}^-\otimes t^{\nu(h_{\alpha-\alpha_{i,s}}
)+\max\{\lambda_1(h_{\alpha})-1,\lambda_2(h_{\alpha})-1\}})w_{\lambda_1+\nu,\lambda_2+\nu}.$$
In the third case if $\lambda(h_{j-1})=1$ or if $\lambda(h_{j-1})=0$ and  there exists $s>j$ minimal with $\lambda(h_s)=1$ then the result follows as before by working with the pairs $(\alpha_{j-1,j}, \alpha-\alpha_{j-1,j})$ and $(\alpha_{j-1,s}, \alpha-\alpha_{j-1,s})$ respectively.
Otherwise we have $\max\{\lambda_1(h_\alpha), \lambda_2(h_\alpha)\} =2$ and the result follows from $$(x_{j-1,n-1}^-\otimes t^{\nu(h_{j-1,n-1})+1})w_{\lambda_1+\nu,\lambda_2+\nu} = 0 = (x_{j,n}^-\otimes t^{\nu(h_{j,n})+1})w_{\lambda_1+\nu,\lambda_2+\nu}.$$
\end{pf}

\subsection{}  To complete the proof of the second step, we now show that \eqref{gdem} also holds. By the preceding Lemma it suffices to  prove that  $$(*)\ \  (x_\alpha^-\otimes t^{\nu(h_\alpha)+\max\{(\lambda_1-\beta_\lambda)(h_\alpha),\lambda_2(h_\alpha)\}}) (x_{\beta_\lambda}^-\otimes t^{(\nu+\lambda_1)(h_{\beta_\lambda})-1})w_{\lambda_1+\nu,\lambda_2+\nu} =0 $$ for roots of the form $\alpha=\alpha_{i,j}$ with $$(**)\ \ (\lambda-\beta_\lambda)(h_i)= 1 =(\lambda-\beta_\lambda)(h_j),\ \ (\lambda-\beta_\lambda)(h_s)=0,\ \ i<s<j,$$ and note that this forces $\beta_\lambda(h_i)\le 0$ and $\beta_\lambda(h_j)\le 0$. We claim that $(\beta_\lambda,\alpha)\le 0$. Otherwise we would have that $\beta_\lambda-\alpha\in R$ and since $\alpha-\beta_\lambda\notin R^+$ we must have $\beta_\lambda-\alpha\in R^+$. It is elementary to see that this is only possible  if  $\alpha_{i,j}=\alpha_{p',r}$ or $\alpha_{p+1,r}$ for some $r\ne p$ and $p'< r\le n$. Since $\beta_{p',p+1}(h_{p'})=1-\delta_{p',p}, $ and $\beta_{p',p+1}(h_{p+1})=1$ we see that $(**)$  forces
$p'=p$ and $\alpha=\alpha_{p,r}$. However, if $p'=p$ then we have $\lambda(h_p)=0$ by definition and $\beta_\lambda(h_p)=0$ which  again contradicts $(**)$ and the claim is proved.  

If $\alpha+\beta_\lambda\notin R^+$, i.e. $(\alpha,\beta_\lambda)=0$ then
$$(x_\alpha^-\otimes t^{\nu(h_\alpha)+\max\{(\lambda_1-\beta_\lambda)(h_\alpha),\lambda_2(h_\alpha)\}})w_{\lambda_1+\nu, \lambda_2+\nu}=0,$$ is a defining relation and $(\ast)$ follows  since $[x_\alpha^-,x_{\beta_\lambda}^-] = 0$. If $(\alpha,\beta_\lambda)<0$ then $\alpha+\beta_\lambda \in R^+$ and  $[x_\alpha^-,x_{\beta_{\lambda}}^-]=x^-_{\alpha+\beta_\lambda}$. It suffices to show that 
$$(x_{\alpha+\beta_\lambda}^-\otimes t^{\nu(h_\alpha)+\max\{(\lambda_1-\beta_\lambda)(h_\alpha),\lambda_2(h_\alpha)\}+(\nu+\lambda_1)(h_{\beta_\lambda})-1})w_{\lambda_1+\nu, \lambda_2+\nu}=0.$$
\noindent But this holds since 
\begin{align*}
    \max\{\lambda_1(h_{\alpha+\beta_\lambda}),\lambda_2(h_{\alpha+\beta_\lambda})\} &=\max\{\lambda_1(h_\alpha) +
    3-\delta_{p',p},\lambda_2(h_\alpha)+1-\delta_{p',p}\}\\ 
    &\le \max\{\lambda_1(h_\alpha) + 3-\delta_{p',p},\lambda_2(h_\alpha)+2-\delta_{p',p}\}\\ &=\max\{\lambda_1(h_\alpha)+1,\lambda_2(h_\alpha)\}+2-\delta_{p',p}.
    \end{align*}
    This completes the proof of Proposition \ref{res}.

 \section{  Connection with the Representation theory of Quantum Affine Algebras}\label{s4}
 In this section we explain very briefly the relationship between our results and the representation theory of quantum affine algebras. We also discuss the connection with the papers \cite{HL1, HL2} on monoidal categorification.
 \subsection{
 }  Let  $\mathbb C(q)$ be   the field of rational functions in an indeterminate $q$ and set  $\bold A=\mathbb Z[q,q^{-1}]$. Let    $\bold U_q(\widehat{\lie g})$ be the  quantized enveloping algebra  (defined over $\mathbb C(q)$)  associated to $\widehat{\lie g}$.    Let   $\bold U_\bold A(\widehat{\lie g})$ be the  $\bold A$-form of   $\bold U_q(\widehat{\lie g})$ and recall that it is  a  free $\bold A$-module and  $$ \bold U_q(\widehat{\lie g})\cong \bold U_\bold A(\widehat{\lie g})\otimes_{\bold A}\mathbb C(q).$$  Regarding  $\mathbb C$ as an $\bold A$-module by letting $q$ act as 1 we have
 that $ \bold U_\bold A(\widehat{\lie g})\otimes_\bold A\mathbb C$ is an  algebra over $\mathbb C$ which has the universal enveloping algebra  $\bold U(\widehat{\lie g})$ as a  canonical quotient.
Finally, recall that  $\bold U_q(\widehat{\lie g})$ is a Hopf algebra and that   $\bold U_\bold A(\widehat{\lie g})$ is a  Hopf subalgebra.

\medskip

Let  $\cal P_{\mathbb Z}^+$ be the free abelian monoid  generated by elements $\{\bomega_{i,q^r}: 1\le i\le n,\  r\in\mathbb Z\}$ and let $\wt:\cal P^+\to P^+$   be the morphism of monoids given by $\wt\bpi=\sum_{i=1}^n(\deg\pi_i) \omega_i$.
\begin{defn}\label{calP} Let $\cal P^+_\mathbb Z(1)$ be the subset of $\cal P^+_\mathbb Z$ containing the identity of the  and elements of the form  $\bomega_{i_1,a_1}\cdots\bomega_{i_k,a_k}$, where $1\le i_1<i_2<\cdots< i_k\le n$,  $a_j\in q^{\mathbb Z}$ for   $1\le k\le n$, and \begin{gather*}\label{defcalp1}a_{j}/a_{{j+1}} =q^{\pm(i_{j+1}-i_j+2)},\ \ k\ge 2,\ \\ a_{j}/a_{{j+1}}= q^{\pm(i_{j+1}-i_j+2)}\implies a_{{j+1}}/a_{{j+2}}= q^{\mp(i_{j+2}-i_{j+1}+2)}\end{gather*} where the second requirement holds for all $j\le k-2$ if $(i_{k-1},i_k)\ne (n-1,n)$  and if $(i_{k-1},i_k)=(n-1,n)$ then it holds for $j\le k-3$ and we require $a_{{k-1}}=a_k$. \hfill\qedsymbol
\end{defn}
Clearly $\wt\cal P^+_\mathbb Z(1)= P^+(1)$. 

\subsection{} To each  element  $\bpi\in \cal P^+_{\mathbb Z}$ one can associate an (unique up to isomorphism) irreducible finite dimensional representation $[\bpi]$ of $\bold U_q(\widehat{\lie g})$.
The trivial representation corresponds to the  identity of the monoid.
Given $\bpi,\bpi'\in\cal P^+_{\mathbb Z}(1)$, the   tensor product $[\bpi]\otimes [\bpi']$ is generically irreducible and isomorphic to $[\bpi \bpi']$. However,  necessary and sufficient conditions for this to hold are not known outside the case of $\lie{sl}_2$ and motivates the interest in understanding the prime irreducible representations.
Recall that a representation is prime if it cannot be written  as a tensor product of two non-trivial representations. 
 The following result is not hard to prove (see for instance \cite{BC, BCM} for a similar statement for $A_n$).
\begin{lem} The module $[\bpi]$ is prime for all $\bpi\in\cal P_\mathbb Z^+(1)$.\hfill\qedsymbol\end{lem}
\noindent The results of
     \cite{Cha01, CPweyl}, show  that the representation $[\bpi]$  admits an $\bold A$-form denoted $[\bpi_\bold A]$ and $[\bpi_\bold A]\otimes_{\bold A} \mathbb C$ is an indecomposable and usually reducible module for  the enveloping algebra $\bold U(\widehat{\lie g})$ and hence also for the current algebra $\lie g[t]$.  Moreover if we pull-back this representation via  the  automorphism of  $\lie g[t]\to\lie g[t]$ sending $x\otimes t^r \to x\otimes (t-1)^r$, $x\in\lie g$, $r\in\mathbb Z_+$ we get a representation $[\bpi_{\mathbb C}]$  of $\lie g[t]$ and   $$(\lie g\otimes t^N\mathbb C[t])[\bpi_{\mathbb C}]=0,\ \ \  N>>0.$$
     \subsection{} We now explain how to deduce the following result. 
     \begin{thm}\label{gradedlimit}
    For $\bpi\in\cal P^+_{\mathbb Z}(1)$, there exists an isomorphism of $\lie g[t]$-modules
     $$[\bpi_{\mathbb C}]\cong D(\lambda_1,\lambda_2),$$ where $\wt\bpi=\lambda$ and $(\lambda_1,\lambda_2)$ is the interlacing pair associated to $\lambda$. In particular $[\bpi_{\mathbb C}]$ acquires the structure of a graded $\lie g[t]$-module.
     \end{thm}
     The theorem  can be  proved by   the same methods as the  ones used in  \cite{BCM} for $\lie{sl}_{n+1}.$
      It uses  the following   idea developed in \cite{moura} where  similar questions were studied for a different family of irreducible representations.
      
Suppose that $\bpi_1,\bpi_2\in\cal P_\mathbb Z^+$ are such that we have an injective map of $\bold U_q(\widehat{\lie g})$-modules $[\bpi_1\bpi_2]\to [\bpi_1]\otimes [\bpi_2]$. Since $\bold U_\bold A(\widehat{\lie g})$ is a Hopf subalgebra of $\bold U_q(\widehat{\lie g})$, we get an injective map $ [(\bpi_1\bpi_2)_\bold A]\to [(\bpi_1)_{\bold A}]\otimes  [(\bpi_2)_{\bold A}]$.  It was shown in \cite[Lemma 2.20, Proposition 3.21]{moura} that tensoring with $\otimes_\bold A\mathbb C$ and pulling back by the automorphism of $\lie g[t]$ induced by  $t\to t-1$,   gives rise to a map of $\lie g[t]$-modules $[\bpi_{\mathbb C}]\to [(\bpi_1)_{\mathbb C}]\otimes [(\bpi_2)_{\mathbb C}]$. 
Here are the main steps in proving the theorem. Retain the notation of the theorem.\\

We   prove that the choice of parameters in Definition \ref{calP} along with Theorem \ref{main} implies that there exists a surjective map of $\lie g[t]$-modules $ D(\lambda_1,\lambda_2)\to [\bpi_{\mathbb C}]\to 0$.\\

Next note that one can define canonically,  elements $\bpi_1,\bpi_2\in\cal P_\mathbb Z^+(1)$  with  $\wt\bpi_s=\lambda_s$, $s=1,2$ and $\bpi=\bpi_1\bpi_2$. Moreover, using the results in  (\cite{CPweyl, FL}), one  has the following isomorphisms of graded $\lie g[t]$-modules, $$[(\bpi_1)_{\mathbb C}]\cong_{\lie g[t]}  D(1,\lambda_1),\qquad\ [(\bpi_2)_{\mathbb C}]\cong_{\lie g[t]} D(1,\lambda_2).$$

The results in \cite{CBraid} show that there is an
 injective map of $\bold U_q(\widehat{\lie g})$-modules $[\bpi]\to [\bpi_1]\otimes [\bpi_2]$. As discussed earlier, one can use \cite[Lemma 2.20, Proposition 3.21]{moura}  to show that the image of the  induced map $$ [\bpi_{\mathbb C}]\to [(\bpi_1)_{\mathbb C}]\otimes [(\bpi_2)_{\mathbb C}]\cong D(1,\lambda_1)\otimes D(1,\lambda_2)$$ is $D(\lambda_1,\lambda_2)$. Hence we have a composition of surjective maps $$D(\lambda_1,\lambda_2)\twoheadrightarrow[\bpi_{\mathbb C}]\twoheadrightarrow D (\lambda_1,\lambda_2).$$  and Theorem \ref{gradedlimit} is proved.

\subsection{The connection with the category $\cal C_\xi$. }\hfill

 We discuss the relationship of our work with that of  \cite{HL1, HL2, nakcluster} and restrict ourselves to $D_n$.
Let  $\xi: \{1,2,\cdots, n\}\to \mathbb Z$ and define the bipartite quiver on the Dynkin diagram of $D_n$ given by  $$\xi(i)=\xi(i+1)\pm 1,\ \ \xi(i)=\xi(i+2),\ \ \xi(n-1)=\xi(n).$$
Let $\cal C_\xi$ be  the full subcategory of  finite-dimensional representations  of $\bold U_q(\widehat{\lie g})$ defined as follows: an object of $\cal C_\xi$   has all its Jordan--Holder components of the form $[\bpi]$, where  $\bpi\in\cal P^{+}$ is a  product of terms of   the form  $\bomega_{i,a}$,  $a\in\{q^{\xi(i)\pm 1}\}$, $1\le i\le n$.  The following was proved  in \cite{HL1,HL2} for  $D_4$ and in \cite{nakcluster} for $D_n$.

\begin{thm}\label{hl}  The category $\cal C_\xi$ is closed under taking tensor products. The Grothendieck ring of $\cal C_\xi$ is a  monoidal categorification of a cluster algebra of type $D_n$ with $n$ frozen variables. \hfill\qedsymbol
\end{thm}
Recall that the cluster variables in a cluster algebra of type  $D_n$ are indexed by elements of  $R^+\cup\{-\alpha_i: 1\le i\le n\}.$ The cluster variable associated to a root $\alpha_{i,j}$ of $D_n$ corresponds to a representation $[\bpi]$ where $\bpi\in\cal P^+_{\mathbb Z}(1)$ is given by $\bpi=\bomega_{i,a_i}\bomega_{i+1,a_{i+1}}\cdots\bomega_{j,a_j}$ and the $a_ks$ are determined by requiring that $\bpi\in\cal P^+_{\mathbb Z}(1)$ and $a_i=q^{\xi(i)\pm 1}$  if $\xi(i)=\xi(i-1)\pm 1$.  Theorem \ref{gradedlimit} gives a character formula for the  prime objects corresponding to these cluster variables. More generally if we relax the condition that $\xi$ define a bipartite quiver, then the  elements of $\cal P^+_{\mathbb Z}(1)$ should correspond to cluster variables in a suitable monoidal categorification.  This was done in detail for $A_n$ in \cite{BC} and it was  non trivial to identify the cluster variable corresponding to $[\bpi]$ when $\xi$ is not bipartite or monotonic.

For $D_n$, even in the bipartite case, it seems to be much more difficult to give the character of the prime representations corresponding to roots of the form $\beta_{i,j}$. However, preliminary calculations suggest that the graded limits of these will also admit a flag where the successive quotients are Demazure modules, but possibly of level bigger than two. We hope to return to these problems elsewhere.

\end{document}